%% file: main.tex
\documentclass[12pt,reqno]{amsart}
\usepackage[utf8]{inputenc} 
\usepackage[T1]{fontenc}    
\usepackage{lmodern}  

\usepackage{graphics,color}
\usepackage{amssymb,amsmath,amsthm,amscd}
\usepackage{latexsym,verbatim,graphicx,amsfonts}
\usepackage{hyperref}
\usepackage[mathscr]{euscript}
\usepackage{amsmath, amsthm, amssymb}
\usepackage{dsfont}
\usepackage{enumitem}
\usepackage{mathrsfs}
\usepackage{mathtools}

\theoremstyle{plain}
\newtheorem{theorem}{Theorem}[section]

\newtheorem{corollary}[theorem]{Corollary}
\newtheorem{lemma}[theorem]{Lemma}
\newtheorem{proposition}[theorem]{Proposition}
\theoremstyle{definition}

\newtheorem{definition}[theorem]{Definition}

\theoremstyle{remark}
\newtheorem{rem}[theorem]{Remark}

\DeclareMathOperator{\spn}{span}

\newcommand{\norm}[1]{\lVert#1\rVert}

\newcommand{\vertiii}[1]{{\left\vert\kern-0.25ex\left\vert\kern-0.25ex\left\vert #1 
    \right\vert\kern-0.25ex\right\vert\kern-0.25ex\right\vert}}

\font\sstext=ecss1000
\font\sssub=ecss1000 at 7pt
\font\sssubsub=ecss1000 at 5pt

\newfam\ssfam
\textfont\ssfam=\sstext
\scriptfont\ssfam=\sssub
\scriptscriptfont\ssfam=\sssubsub



\begin{document}

\title[Operators on injective tensor products]{Operators on injective tensor products of separable Banach spaces and spaces with few operators}
\author[A. Acuaviva]{Antonio Acuaviva}
\address{School of Mathematical Sciences,
Fylde College,
Lancaster University,
LA1 4YF,
United Kingdom} \email{ahacua@gmail.com}

\date{\today}

\begin{abstract}
We give a characterization of the operators on the injective tensor product $E \hat{\otimes}_\varepsilon X$ for any separable Banach space $E$ and any (non-separable) Banach space $X$ with few operators, in the sense that any operator $T: X \to X$ takes the form $T = \lambda I + S$ for a scalar $\lambda \in \mathbb{K}$ and an operator $S$ with separable range. This is used to give a classification of the complemented subspaces and closed operator ideals of spaces of the form $C_0(\omega \times K_\mathcal{A})$, where $K_\mathcal{A}$ is a locally compact Hausdorff space induced by an almost disjoint family $\mathcal{A}$ such that $C_0(K_\mathcal{A})$ has few operators.
\end{abstract}

\maketitle

\bigskip
\input{1-Introduction_and_results}
\bigskip
\input{2-Organization_and_Notation}

\bigskip
\input{3-Proof_Main_Theorem}
\bigskip
\input{4-Applications}
\bigskip

\noindent\textbf{Note added in proof.} After presenting our results in a seminar at Lancaster University, Dr.~Matthew Daws extended Theorem~\ref{th: main-tensor-prod} to operators $T: E \hat{\otimes}_\alpha X \to F \hat{\otimes}_\alpha X$, where $\alpha(\cdot)$ is any uniform crossnorm and $F$ has the bounded approximation property. The core of the argument involves replacing the evaluation and inclusion operators $\pi_k$ and $\iota_e$ with \emph{slice maps}, and employing the principle of local reflexivity to obtain an analogue of Lemma~\ref{lmm: approx2}. \bigskip

\noindent\textbf{Acknowledgements.} This paper forms part of the author’s PhD research at Lancaster University under the supervision of Professor N. J. Laustsen. The author extends sincere gratitude to Professor Laustsen for his insightful comments and valuable suggestions, which greatly improved the manuscript's presentation. Furthermore, his identification of Lemma \ref{lmm: factorization-projections} significantly simplified the proof of Theorem \ref{th: complemented-subspaces}. We also thank the anonymous referee for helpful comments that enhanced the clarity of the document.

He also acknowledges with thanks the funding from the EPSRC (grant number EP/W524438/1) that has supported his studies.


\input{bibliography}
\end{document}

%% file: 1-Introduction_and_results.tex
\section{Introduction}

A common characteristic of many Banach spaces is the rich structure of their spaces of operators. This has, in turn, led to the question of how simple the space of operators
between two Banach spaces can be. Efforts to explore this question have spurred a rich vein of research, culminating in the construction of several exotic Banach spaces.

Among the most notable examples are the space constructed by Gowers and Maurey \cite{gowers1993unconditional}, that does not contain any unconditional basic sequence, and the space built by Argyros and Haydon, which solved the scalar-plus-compact problem \cite{argyros2011hereditarily}. \\

In the non-separable setting, a similar question has been explored, focusing on spaces $X$ where every operator is a separable perturbation of a scalar multiple of the identity; that is, $T = \lambda I + S$, where $S$ is an operator with separable range. Such spaces, which we will refer to as \emph{spaces with few operators}, were first constructed by Shelah \cite{shelah1978banach} under the assumption of the diamond axiom, which was later removed by Shelah and Stepr\={a}ns \cite{shelah1988banach}. Another space of this form, realized as a dual space, was constructed by Argyros and Tolias \cite{argyros2004methods}, while a reflexive example was built by Wark \cite{wark2001non} (note that this implies that its dual also has few operators). 

This question has recently attracted significant interest in the context of spaces of continuous functions, leading to the construction of many exotic locally compact spaces. In this vein, Koszmider constructed, under the Continuum Hypothesis or Martin's Axiom, a scattered, locally compact Hausdorff space $L$ such that $C_0(L)$ admits few operators \cite{koszmider2005decompositions}. These set-theoretic assumptions were later removed by Koszmider and Laustsen \cite{koszmider2021banach}. Additionally, under the assumption of Ostaszewski’s $\clubsuit$-principle, Koszmider and Zieliński \cite{koszmider2011complementation} constructed a different example of a $C_0(L)$ spaces with few operators. This construction was further developed in work of Candido \cite{candido2021banach} and \cite{candido2024few}. \\

Given the simple structure of the operators on these spaces, a natural question arises: can they be used to construct new spaces that still retain easily characterizable operators? We show that this is the case when considering the injective tensor product with a separable Banach space $E$. This result becomes particularly interesting when $X = C_0(L)$ and $E = C(M)$, where $M$ is a compact metric space and $L$ is a locally compact Hausdorff space since in this case $C(M) \hat{\otimes}_\varepsilon C_0(L) \cong C_0(M \times L)$.

We present our main result in slightly more general terms, as this does not complicate the proof. In the following theorem, $\hat{\otimes}_\alpha$ denotes any tensor completion for a \emph{uniform crossnorm} $\alpha(\cdot)$, see Definition \ref{def: uniform-cross-norms}.

\begin{theorem}\label{th: main-tensor-prod}
    Let $E, F$ be separable Banach spaces and $X$ a non-separable Banach space with few operators. Then any operator $T: E \hat{\otimes}_\alpha X \to F \hat{\otimes}_\varepsilon X$ can be uniquely expressed as $T = U \hat{\otimes} I + S$, where $U: E \to F$ and $S$ has separable range.
\end{theorem}

The essence of the previous theorem is as follows: when we take the injective tensor product of a separable Banach space with a space with few operators, we obtain, up to separable elements, the original space we started with. In particular, when taking the injective tensor product of separable Banach spaces with spaces with few operators, the result is an injective operation, in the following sense.

\begin{corollary}\label{cor: inj-operation}
    Let $E, F$ be separable Banach spaces and $X$ be a non-separable Banach space with few operators. Then
    \begin{equation*}
        E \hat{\otimes}_\varepsilon X \sim F \hat{\otimes}_\varepsilon X \iff E \sim F.
    \end{equation*}
\end{corollary}

As previously mentioned, we are particularly interested in the case of spaces of continuous functions. Theorem \ref{th: main-tensor-prod} becomes especially relevant when we take $E = C(M_1)$, $F = C(M_2)$, $\hat{\otimes}_\alpha = \hat{\otimes}_\varepsilon$ the injective tensor product and $X = C_0(L)$ for some locally compact Hausdorff space $L$. This gives information about operators on spaces of the form $C_0(M \times L)$.

\begin{corollary}\label{corollary-operators-CK}
    Let $L$ be a non-metrizable locally compact Hausdorff space such that $C_0(L)$ admits few operators, and $M_1, M_2$ be compact metric spaces. Then any operator $T: C_0(M_1 \times L) \to C_0(M_2 \times L)$ can be uniquely expressed as
    \begin{equation*}
        T = U \hat{\otimes} I + S,
    \end{equation*}
    where $U: C(M_1) \to C(M_2)$ and $S$ has separable range.
\end{corollary}

The previous corollary generalizes Candido's result \cite[Theorem 1.1]{candido2021banach}, proved under the assumption of Ostaszewski’s $\clubsuit$-principle. It extends the result to also include the case $M = [0,1]$ and demonstrates that the key property is that $C_0(L)$ has few operators. For example, it applies to the space built by Koszmider and Laustsen, whose construction does not require additional set-theoretic assumptions. \\

As an application of Theorem \ref{th: main-tensor-prod}, we expand the list of spaces of continuous functions with completely understood complemented subspaces, obtaining the following.

\begin{theorem}\label{th: complemented-subspaces}
    Let $K_\mathcal{A}$ be a locally compact Hausdorff space induced by an almost disjoint family $\mathcal{A} \subseteq [\mathbb{N}]^\omega$ of cardinality the continuum such that $C_0(K_\mathcal{A})$ admits few operators.  
    Then any infinite-dimensional complemented subspace of $C_0(\omega \times K_\mathcal{A})$ is isomorphic to exactly one of the following: $c_0$, $C_0(K_\mathcal{A})^n$ for some $n \in \mathbb{N}$, or $C_0(\omega \times K_\mathcal{A})$.  

\end{theorem}

\begin{rem}
    Note that, if $C(L)$ is a Grothendieck space such that its complemented subspaces are classified, then we also obtain a classification of the complemented subspaces of $C_0((\omega \times K_\mathcal{A}) \sqcup L) \cong C_0((\omega \times K_\mathcal{A})) \oplus C(L) $, see \cite{johnson2016closed}. For example, we get a classification of the complemented subspaces of $\ell^c_\infty(\Gamma) \oplus C(\omega \times K_\mathcal{A})$ for any set $\Gamma$.
\end{rem}

We believe that it is possible to extend the previous classification of complemented subspaces to $C_0(\omega^\omega \times K_\mathcal{A})$, by using arguments similar to those of Benyamini's classification of the complemented subspaces of $C(\omega^\omega)$ \cite{benyamini1978extension}. However, the proof becomes significantly technical and involved, and we are not aware of any immediate applications. We conjecture that it may even be possible to extend the classification to $C_0(M \times K_\mathcal{A})$ for $M$ a compact metric space, in terms of the complemented subspaces of $C(M)$, as in \cite[Theorem 1.3]{candido2021banach}, though we are unable to determine a general approach to achieve this. \\

Finally, we can fully describe the lattice of closed operator ideals of $\mathscr{B}(C_0(\omega \times K_\mathcal{A}))$. Throughout, by an ideal we mean a two-sided ideal in the operator algebra.

For a Banach space $Y$, we denote by $\mathscr{K}(Y)$ the ideal of compact operators, $\mathscr{X}(Y)$ denotes the ideal of operators with separable range, and if $Z$ is another Banach space, we define the ideal $\mathscr{G}_Z(Y)$ by
\begin{equation*}
    \mathscr{G}_Z(Y) = \spn \{ ST: T \in \mathscr{B}(Y, Z), S \in \mathscr{B}(Z, Y)\},
\end{equation*}
where $\mathscr{B}(Y, Z)$ and $\mathscr{B}(Z, Y)$ denote the space of all operators from $Y$ to $Z$ and $Z$ to $Y$ respectively.

\begin{theorem}\label{th: ideal-lattice}
    Let $K_\mathcal{A}$ be a locally compact Hausdorff space induced by an almost disjoint family $\mathcal{A} \subseteq [\mathbb{N}]^\omega$ of cardinality the continuum such that $C_0(K_\mathcal{A})$ admits few operators. Then the lattice of closed ideals of $\mathscr{B}(C_0(\omega \times K_\mathcal{A}))$ is given by
    \begin{align*}
        \{0\} &\subsetneq \mathscr{K}(C_0(\omega \times K_\mathcal{A})) \subsetneq \mathscr{G}_{c_0}(C_0(\omega \times K_\mathcal{A})) = \mathscr{X}(C_0(\omega \times K_\mathcal{A}))\\
        &\subsetneq \overline{\mathscr{G}}_{C_0(K_\mathcal{A})}(C_0(\omega \times K_\mathcal{A})) \subsetneq \mathscr{B}(C_0(\omega \times K_\mathcal{A})).
    \end{align*}
\end{theorem}

\begin{rem} 
    It is worth noting that even more can be deduced from this result. Since the operators on $c_0$ are particularly well-behaved, the same holds for the operators on $C_0(\omega \times K_{\mathcal{A}})$. Specifically, from Theorem~\ref{th: main-tensor-prod}, the Banach algebras 
    \begin{equation*} \mathscr{B}(C_0(\omega \times K_{\mathcal{A}})) / \mathscr{X}(C_0(\omega \times K_{\mathcal{A}})) \quad \text{and} \quad \mathscr{B}(c_0) \end{equation*} 
    are isometrically isomorphic. This allows us to conclude that many properties of operators on $c_0$ also hold for those on $C_0(\omega \times K_{\mathcal{A}})$. For example, it follows that the quotient \begin{equation*} \mathscr{B}(C_0(\omega \times K_{\mathcal{A}})) / \mathscr{I} \end{equation*} has a unique algebra norm for every closed ideal $\mathscr{I}$ of $\mathscr{B}(C_0(\omega \times K_{\mathcal{A}}))$, in the spirit studied in \cite{arnott2023uniqueness}.
\end{rem}
\begin{rem}
    Note that $C_0(\omega \times K_\mathcal{A}) \cong C_0(K_\mathcal{B})$ where $K_\mathcal{B}$ is a locally compact Hausdorff space induced by the almost disjoint family $\mathcal{B}$, obtained by breaking up $\mathbb{N}$ into countably many copies of itself and using the order isomorphism to create a copy of $\mathcal{A}$ in each of these copies of $\mathbb{N}$. Thus, the previous argument shows that the unital Banach algebra $\mathscr{B}(c_0)$ can be realized as a quotient $\mathscr{B}(C_0(K_\mathcal{B}))/\mathscr{X}(C_0(K_\mathcal{B}))$, addressing a question of Koszmider and Laustsen \cite[Question 46 (iii)]{koszmider2021banach}.
\end{rem}

%% file: 2-Organization_and_Notation.tex
\section{Organization and notation}

All normed spaces and algebras are over the scalar field $\mathbb{K}$, either the real or complex numbers, and we adhere to standard notational conventions. The term \emph{operator} will refer to a bounded linear map between normed spaces. For two Banach spaces $X$ and $Y$, we write $X \cong Y$ to mean that $X$ and $Y$ are isometrically isomorphic, and write $X \sim Y$ to mean that they are isomorphic.  We denote by $B_X$ the unit ball of a Banach space $X$. 

More specialized notation will be introduced as and when needed. \\

We give a brief layout of the structure of this paper. In subsection \ref{sub-tensor-products}, we introduce notation regarding tensor products as well as basic results that will be necessary for the proof of Theorem \ref{th: main-tensor-prod}.

In Section \ref{sec: proof-Theorem}, we present the proof of Theorem \ref{th: main-tensor-prod}. First, we prove the result in the case $F = C(M)$. For this, the key is to construct the operator $U: E \to C(M)$ and then verify that $T - U \hat{\otimes} I$ has separable range. The general case when $F$ is any separable Banach space will readily follow using the universality of $C[0,1]$ for separable Banach spaces, together with the fact that the injective tensor product preserves subspaces. We also provide a proof of Corollary \ref{cor: inj-operation}.

In Section \ref{sec: proof-applications}, we explore some applications of our theorem in the context of Banach spaces induced by an almost disjoint family, including the proofs of Theorems \ref{th: complemented-subspaces} and \ref{th: ideal-lattice}.

\subsection{Overview of Tensor Products}\label{sub-tensor-products}

We briefly review some standard facts about tensor products of Banach spaces in general, and the injective tensor product in particular, mainly following \cite{ryan2002introduction}. We assume familiarity with the algebraic tensor product $X \otimes Y$ of two vector spaces $X$ and $Y$, and refer to \cite[Chapter~1]{ryan2002introduction} for the basic definitions. Recall that if $E$, $F$, $X$, and $Y$ are vector spaces and $U: E \to F$, $R: X \to Y$ are linear maps, then we can define the tensor product linear map $U \otimes R: E \otimes X \to F \otimes Y$ which satisfies
\begin{equation*}
     (U \otimes R)(e \otimes x) = (U e) \otimes (R x)
\end{equation*}
for every $e \in E$ and $x \in X$. 

Given two Banach spaces $X$ and $Y$, their algebraic tensor product $X \otimes Y$ is a vector space that can be equipped with various norms. The theory of tensor products of Banach spaces studies how to define such norms and the properties of the resulting Banach spaces obtained by completing with respect to these norms.

If $\alpha(\cdot)$ is a norm on the algebraic tensor product $X \otimes Y$, then $X \hat{\otimes}_{\alpha} Y$ is defined as the completion of $X \otimes Y$ under the norm $\alpha$. Among the most important norms are the projective norm $\pi(\cdot)$ (see \cite[Chapter~2]{ryan2002introduction}) and the injective norm $\varepsilon(\cdot)$ (see \cite[Chapter~3]{ryan2002introduction}). We recall that for an element $u = \sum_{i=1}^n x_i \otimes y_i \in X \otimes Y$, we have
\begin{equation*}
    \varepsilon(u) = \sup \left \{\left|\sum_{i=1}^n \varphi(x_i)\psi(y_i) \right|: \varphi \in B_{X^*}, \psi \in B_{Y^*} \right\},
\end{equation*}
and that this norm is well-defined and independent of the representation we have chosen for $u$. We also have the following elementary property: $\varepsilon(x \otimes y) = \norm{x} \norm{y}$ for every $x \in X$ and $y \in Y$ (see \cite[Proposition 3.1]{ryan2002introduction}). An important fact that we will use is that the injective tensor product preserves subspaces. We also have the following, see \cite[Proposition 3.2]{ryan2002introduction}.

\begin{proposition}\label{prop: injective-preserve-subspaces}
    Let $U: E \to F$ and $R: X \to Y$ be operators. Then there exists a unique operator $U \hat{\otimes}R: E \hat{\otimes}_\varepsilon X \to F \hat{\otimes}_\varepsilon Y$ such that $(U \hat{\otimes}R)(e \otimes x) = (Ue) \otimes (Rx)$ for every $e \in E$ and $x \in X$. Furthermore $\norm{U \hat{\otimes}R} = \norm{U} \norm{R}$.
\end{proposition}

We note that, in the previous proposition, the operator $U \hat{\otimes}R$ is obtained as the continuous extension of the linear map $U \otimes R: E \otimes X \to F \otimes Y$. We also recall the following well-known result, see \cite[Section 3.2]{ryan2002introduction}.

\begin{proposition}
    Let $K$ and $L$ be compact Hausdorff spaces and $X$ be a Banach space. Then
    \begin{equation*}
        C(K) \hat{\otimes}_{\varepsilon} X \cong C(K, X),
    \end{equation*}
    where the identification for elementary tensors is given by 
    \begin{equation*}
        (f \otimes x)(k) = f(k) x \quad \text{for all } f \in C(K),\ x \in X, \text{ and } k \in K.
    \end{equation*}
    Similarly,
    \begin{equation*}
        C(K) \hat{\otimes}_{\varepsilon} C(L) \cong C(K \times L).
    \end{equation*}
\end{proposition}

We have seen that for the injective tensor product, given operators $U: E \to F$ and $R: X \to Y$, we can define an operator
\begin{equation*}
    U \hat{\otimes} R : E \hat{\otimes}_\varepsilon X \to F \hat{\otimes}_\varepsilon Y,
\end{equation*}
by extending the linear map $U \otimes R$.
For our purposes, it is important to ask for which norms $\alpha(\cdot)$ on the algebraic tensor product the linear map $U \otimes R$ extends to an operator
\begin{equation*}
    U \hat{\otimes} R : E \hat{\otimes}_\alpha X \to F \hat{\otimes}_\varepsilon Y,
\end{equation*}
where we emphasize that we are mapping into the injective tensor product $F \hat{\otimes}_\varepsilon Y$ from the $\alpha$-tensor product $E \hat{\otimes}_\alpha X$.

In this paper, we restrict our attention to a particular well-behaved family of norms (see \cite[Section~6.1]{ryan2002introduction}), which in particular guarantees the existence of such extensions.

\begin{definition}\label{def: reasonable-cross-norms}
    Let $X$ and $Y$ be Banach spaces. A norm $\alpha(\cdot)$ on $X \otimes Y$ is a \emph{reasonable crossnorm} if it has the following properties:
    \begin{enumerate}
        \item $\alpha(x \otimes y) \leq \norm{x}\norm{y}$ for every $x \in X$, $y \in Y$.
        \item For every $\varphi \in X^*$ and $\psi \in Y^*$, the linear functional $\varphi \otimes \psi$ on $X \otimes Y$ is bounded, and $\norm{\varphi \otimes \psi} \leq \norm{\varphi} \norm{\psi}$.
    \end{enumerate}
\end{definition}

\begin{definition}\label{def: uniform-cross-norms}
    A \emph{uniform crossnorm} is an assignment to each pair of Banach spaces $(E,X)$ of a reasonable crossnorm $\alpha(\cdot)$ on $E \otimes X$ such that, for all operators $U: E \to F$ and $R: X \to Y$, we have
    \begin{equation*}
        \norm{U \otimes R} \leq \norm{U} \, \norm{R},
    \end{equation*}
    where the operator norm on the left is computed with respect to the corresponding crossnorm assignments.
\end{definition}

It is well-known that $\alpha(\cdot)$ is a reasonable crossnorm if and only if $\varepsilon(u) \leq \alpha(u) \leq \pi(u)$ for every $u \in X \otimes Y$, see \cite[Proposition 6.1]{ryan2002introduction}. It follows that if $\alpha(\cdot)$ is a uniform crossnorm and $U: E \to F$ and $R: X \to Y$ are operators, we can define the product operator
\begin{equation*}
    U \hat{\otimes} R: E \hat{\otimes}_{\alpha} X \to F \hat{\otimes}_{\varepsilon} Y,
\end{equation*}
so that $\norm{U \hat{\otimes} R} \leq \norm{U} \norm{R}$ and $(U \hat{\otimes} R)(e \otimes x)  = (U \otimes R)(e \otimes x) = (Ue) \otimes (Rx)$ for every $e \in E$ and $x \in X$. Lastly, we note that the injective norm $\varepsilon(\cdot)$ is a uniform cross norm. This is everything we will need from the theory of tensor products in Banach spaces.

%% file: 3-Proof_Main_Theorem.tex
\section{Proof of Theorem \ref{th: main-tensor-prod}}\label{sec: proof-Theorem}

Throughout this section, $X$ will always denote a non-separable Banach space with few operators, meaning that every operator $T: X \to X$ can be expressed as $T = \lambda I + S$, where $\lambda \in \mathbb{K}$ and $S$ is an operator with separable range. Meanwhile, $E$ and $F$ will always refer to separable Banach spaces, while $M$ will be a metrizable compact space.

The strategy for the proof is to exploit the simple structure of operators on $X$. We first assume that $F = C(M)$ and we need to transform an operator $T: E \hat{\otimes}_\alpha X \to C(M, X)$ into a collection of operators from $X$ to $X$, which encapsulate all the information of the original operator. We will need the following definition.

\begin{definition}
    For each $k \in M$, define the \emph{evaluation operator}
    \begin{equation*}
        \pi_k: C(M, X) \to X, \hspace{10pt} g \mapsto g(k).
    \end{equation*}
    For each $e \in E$, define the \emph{inclusion operator}
\begin{equation*}
    \iota_e: X \to E \hat{\otimes}_\alpha X, \hspace{10pt}  x \mapsto e \otimes x.
\end{equation*}
 Note that $\norm{\iota_e} \leq \norm{e}$ and $\norm{\pi_k} = 1$.
\end{definition}

Using the evaluation and inclusion operators, it is easy to obtain operators from $X$ to $X$ using our original operator $T$. We will then use these operators to build $U: E \to C(M)$. 

\begin{definition}
    Let $e \in E$, $k \in M$ and $T: E \hat{\otimes}_\alpha X \to C(M, X)$ be an operator. Then the \emph{$(e,k)-$coordinate} of $T$ is given by
 \begin{equation*}
     T_{e,k} = \pi_kT\iota_e: X \to X,
 \end{equation*}
 so that in particular there is a unique decomposition $T_{e,k} = \lambda_{e,k}I + S_{e,k}$ where $\lambda_{e,k} \in \mathbb{K}$ and $S_{e,k}: X \to X$ is an operator with separable range.
\end{definition}

From now on, fix a countable dense subset $M_0 \subseteq M$. We emphasise that this set will be employed repeatedly throughout this section.

\begin{definition}\label{def: E_0generated}
Let $E_0$ be a countable subset of $E$. We denote by $Y_{E_0}$ the (separable) subspace
\begin{equation*}
        Y_{E_0} := \overline{\operatorname{span}} \bigcup \{ S_{e,k}[X] : e \in E_0,\, k \in M_0 \}.
\end{equation*}
\end{definition}

We have the following elementary observation.

\begin{lemma}\label{lmm: approx-coordinatewise}
    Let $E_0 \subseteq E$ be a countable subset and $Y_{E_0}$ be as in Definition \ref{def: E_0generated}. Then there exist $x \in B_{X}$ and $\varphi \in B_{X^*}$ such that $\varphi(x) \geq 1/2$ and $\varphi(Y_{E_0}) = 0$.

    In particular, for any finite collection of vectors $e_1, \dots, e_n \in E$, there exists $x \in B_X$ such that
    \begin{equation*}
        \left|\sum_{i=1}^n c_i \lambda_{e_i,k_i} \right| \leq 2\left \lVert \sum_{i=1}^n c_i \pi_{k_i} T (e_i \otimes x) \right \lVert
    \end{equation*}
    for any choice $k_1, \dots, k_n \in M_0$ and scalars $c_1, \dots, c_n \in \mathbb{K}$.
\end{lemma}
\begin{proof}
    Note that $Y_{E_0} \subsetneq X$, since the latter is non-separable. Thus an application of the Hahn-Banach theorem gives the first part.

    For the second part, take $E_0 = \{e_1, \dots, e_n\}$. There exists $x \in B_X$ and $\varphi \in B_{X^*}$ such that $\varphi(Y_{E_0}) = 0$ and $\varphi(x) \geq 1/2$. It follows that for any $k_1, \dots, k_n \in M_0$ and $c_1, \dots, c_n \in \mathbb{K}$, we have
    \begin{align*}
         \left|\sum_{i=1}^n c_i \lambda_{e_i,k_i} \right|&\leq  2|\varphi(x)|\left|\sum_{i=1}^n c_i \lambda_{e_i,k_i}\right| = 2\left|\varphi\left(\sum_{i=1}^n c_i \lambda_{e_i,k_i}x + c_i S_{e_i, k_i}x \right) \right|\\
         &\leq  2\norm{\sum_{i=1}^n c_i ( \lambda_{e_i,k_i}x + S_{e_i, k_i}x)} = 2\norm{\sum_{i=1}^n c_i \pi_{k_i} T (e_i \otimes x)}. \qedhere
    \end{align*}
\end{proof}

We can now build the operator $U: E \to C(M)$, using $(\lambda_{e,k})_{e \in E, k \in M_0}$ together with an extension via density.

\begin{lemma}\label{lmm: approx1}
    For each $e \in E$ the function $U_0e: M_0 \to \mathbb{K}$ defined by $U_0e(k) = \lambda_{e,k}$ is uniformly continuous and thus admits a unique continuous extension to $M$, which we denote by $Ue \in C(M)$.
\end{lemma}
\begin{proof}
    Let $e \in E$ be fixed, applying the second part of Lemma \ref{lmm: approx-coordinatewise} for any $k_1, k_2 \in M_0$ we have
    \begin{align*}
        |U_0e(k_1) - U_0e(k_2)| &= |\lambda_{e, k_1} - \lambda_{e, k_2}| \leq 2\norm{\pi_{k_1} T (e \otimes x) - \pi_{k_2} T (e \otimes x)} \\
        &=2\norm{(T (e \otimes x))(k_1) -  (T (e \otimes x))(k_2)},
    \end{align*}  
    and since $T(e \otimes x) \in C(M, X)$ is uniformly continuous by the Heine-Cantor theorem, so is $U_0 e$.
\end{proof}

\begin{lemma}\label{lmm: approx2}
    The map $U: E \to C(M)$, $e \mapsto Ue$, where $Ue$ is given as in Lemma \ref{lmm: approx1}, defines an operator of norm at most $2\norm{T}$.
\end{lemma}
\begin{proof}
    We start by showing that $U$ is bounded, by arguing as in the proof of Lemma \ref{lmm: approx1}. For a fixed $e \in E$, an application of Lemma \ref{lmm: approx-coordinatewise} gives
    \begin{equation*}
        \norm{Ue} = \sup_{k \in M_0} \norm{U_0e(k)} = \sup_{k \in M_0} |\lambda_{e,k}| \leq \sup_{k \in M_0}\norm{\pi_k T (e \otimes x)}\leq 2\norm{T}\norm{e}.
    \end{equation*}

    To show linearity, let $e_1, e_2 \in E$ and $c \in \mathbb{K}$. An application of Lemma \ref{lmm: approx-coordinatewise}, together with the linearity of $T$ and the tensor product, gives $|(\lambda_{e_1, k} + c\lambda_{e_2, k}) - \lambda_{e_1 + ce_2, k}| = 0$ for each $k \in M_0$. Therefore, the continuous function $(Ue_1 + cUe_2) - U(e_1 + ce_2)$ vanishes on the dense set $M_0$, and thus it is the null function. This proves the linearity of $U$.
\end{proof}

Once we have found the operator $U: E \to C(M)$, we simply need to check that it gives the desired representation, that is, we need to check that the operator $S = T - U \hat{\otimes} I$ has separable range. For that, we need the following lemma.

\begin{lemma}\label{lmm: separable}
    Let $Y$ be a subspace of $C(M, X)$. Then $Y$ is separable if and only if $\pi_{k} [Y]$ is separable for each $k \in M_0$.
\end{lemma}
\begin{proof}
    If $Y$ is separable, clearly $\pi_k [Y]$ is separable, so we focus on the reverse implication. Let $Z = \overline{\spn} \bigcup_{k \in M_0} \pi_k [Y]$ which is separable, since by assumption each $\pi_k [Y]$ is. 
    We naturally view $C(M, Z)$ as a closed subspace of $C(M, X)$. Observe that $C(M, Z)$ is also separable, since $C(M, Z) \cong C(M) \hat{\otimes}_{\varepsilon} Z$ and the tensor product of separable spaces is itself separable. Therefore, it is enough to show that $Y \subseteq C(M, Z)$, 

    For this, we need to show that every $y \in Y$ takes values in $Z$. For any $k \in M$, take $k_n \in M_0$, $k_n \to k$. By continuity of $y$, $\pi_{k_n}(y) \to \pi_{k}(y)$. Since $\pi_{k_n}(y) \in Z$ for every $n \in \mathbb{N}$ and $Z$ is a closed subspace, then $\pi_{k}(y) \in Z$, finishing the proof.
\end{proof}

We are ready to prove Theorem \ref{th: main-tensor-prod} in the case $F = C(M)$.

\begin{proof}[Proof of Theorem \ref{th: main-tensor-prod}, for $F = C(M)$]
    Let $U: E \to C(M)$ be the operator given by Lemma \ref{lmm: approx2}. We claim that the operator $S = T - U \hat{\otimes} I$ has separable range. The uniqueness of the representation is clear, once its existence has been shown.

    Thus, we only need to prove that $S$ has separable range. Let $E_0$ be a countable dense set of $E$ and recall that $Y_{E_0}$, denotes the closed linear span of $\bigcup \{S_{e,k}[X]: e \in E_0, k \in M_0\}$. By Lemma \ref{lmm: separable}, we need to show that $\pi_k S[E \hat{\otimes}_\alpha X]$ is separable for each $k \in M_0$, and for this, it is enough to show $\pi_k S[E \hat{\otimes}_\alpha X] \subseteq Y_{E_0}$. 

    Since $E_0$ is dense in $E$, $\spn\{ e \otimes x: e \in E_0, x \in X\}$ is dense in $E \hat{\otimes}_\alpha X $. Therefore, by continuity and linearity, it is enough to show that $\pi_kS(e \otimes x) \in Y_{E_0}$ for each $e \in E_0$ and $x \in X$. We have
    \begin{align*}
        \pi_k S(e \otimes x) &= \pi_k (T - U \hat{\otimes} I)(e \otimes x) = \pi_k T (e \otimes x) - \pi_k((Ue) \otimes x) \\
        &= (\lambda_{e,k}x + S_{e,k}x) - \lambda_{e,k}x = S_{e,k}x \in Y_{E_0},
    \end{align*}
    as required.
\end{proof}
\begin{proof}[Proof of Theorem \ref{th: main-tensor-prod}, general case.]

Take $M = [0,1]$ and let $j: F \hookrightarrow C(M)$ be an isometric embedding. Since the injective tensor product preserves subspaces, we can identify $F \hat{\otimes}_\varepsilon X$ with a closed subspace of $C(M, X) \cong C(M) \hat{\otimes}_\varepsilon X$, specifically $j \hat{\otimes} I: F \hat{\otimes}_\varepsilon X \hookrightarrow C(M) \hat{\otimes}_\varepsilon X$ is an isometric embedding.

Let $T: E \hat{\otimes}_\alpha X \to F \hat{\otimes}_\varepsilon X$ be an operator. By the previous case applied to $(j \hat{\otimes} I) \circ T$, there exists $U: E \to C(M)$ and an operator $S$ with separable range such that
\begin{equation*}
    (j \hat{\otimes} I) \circ T = U \hat{\otimes} I + S.
\end{equation*}
Assume we have shown that $U[E] \subseteq j[F] \subseteq C(M)$, so that $U = j \circ V$ for some $V: E \to F$. Therefore
\begin{equation*}
    U \hat{\otimes} I = (j \circ V) \hat{\otimes} I = (j \hat{\otimes} I) \circ (V \hat{\otimes} I),
\end{equation*}
which combined with the previous equation gives
\begin{equation*}
    (j \hat{\otimes} I) \circ (T - V \hat{\otimes} I) = S.
\end{equation*}
Since $j \hat{\otimes} I$ is an isometry, we obtain that $T - V \hat{\otimes} I$ has separable range, and the result follows. To finish, we only need to show that $U[E] \subseteq j[F]$. 

Since $j[F]$ is closed, it is enough to show that for any $\varepsilon > 0$ and $e \in E$ we have $d(Ue, j[F]) < \varepsilon$. From now on, fix $\varepsilon > 0$ and $e \in E$. Apply the first part of Lemma \ref{lmm: approx-coordinatewise} to $E_0 = \{e\}$, to obtain $x \in B_X$ and $\varphi \in B_{X^*}$ such that $\varphi(x) \geq 1/2$ and $\varphi(Y_{E_0}) = 0$. 

Since $\spn{ \{f \otimes x: f \in F, x \in X \} }$ is dense in $F \hat{\otimes}_\varepsilon X$ we can choose \break $f_1, \dots, f_n \in F$, $x_1, \dots, x_n \in X$ such that
\begin{equation*}
     \norm{T(e \otimes x) - \sum_{i=1}^n f_i \otimes x_i} < \varepsilon/2.  
\end{equation*}
On the other hand, since $j \hat{\otimes} I$ is an isometry and $M_0$ is a dense set in $M=[0,1]$, we can compute this norm as
\begin{align*}
    \sup_{k \in M_0} \norm{\pi_k ((j \hat{\otimes} I)T)(e \otimes x) &- \sum_{i=1}^n j(f_i)(k)x_i} \\
    &= \sup_{k \in M_0} \norm{(\lambda_{e,k}x + S_{e,k}x) - \sum_{i=1}^n  j(f_i)(k)x_i} \\
    & \geq \sup_{k \in M_0} |\varphi((\lambda_{e,k}x + S_{e,k}x) - \sum_{i=1}^n  j(f_i)(k)x_i)| \\
     &= \sup_{k \in M_0} |\varphi(x)\lambda_{e,k} - \sum_{i=1}^n  \varphi(x_i) j(f_i)(k)| \\
     & \geq \frac{1}{2} \sup_{k \in M_0} |\lambda_{e,k} - \sum_{i=1}^n  \frac{\varphi(x_i)}{\varphi(x)} j(f_i)(k)| \\
     &\geq \frac{1}{2}d(Ue, j(F)),
\end{align*}
where the last inequality follows again since $M_0$ is a dense subset of $M = [0,1]$. It follows that $d(Ue, j(F)) < \varepsilon$, as desired.
\end{proof}
Finally, we deduce Corollary \ref{cor: inj-operation}.
\begin{proof}[Proof of Corollary \ref{cor: inj-operation}.]
    Clearly, $E \sim F$ implies $E \hat{\otimes}_\varepsilon X \sim F \hat{\otimes}_\varepsilon X$, so we focus on the reverse implication.
    
    Let $T: E \hat{\otimes}_\varepsilon X \to F \hat{\otimes}_\varepsilon X$ be an isomorphism. Using Theorem \ref{th: main-tensor-prod} we can express $T = U \hat{\otimes} I + S$ where $U: E \to F$ and $S$ has separable range. We claim that $U$ is an isomorphism, which will finish the proof. To establish this, we will show that $U$ is bounded below and that its range $U[E]$ is dense in $F$. These two properties suffice because a bounded below operator is injective with closed range, and if its range is also dense, then $U$ is surjective.

    We start by showing that $U$ is bounded below, so fix $e \in E$ with $\norm{e} = 1$. Since $X$ is non-separable, by transfinite induction we can find $(x_\alpha)_{\alpha < \omega_1} \subseteq B_X$ such that $\norm{x_\alpha - x_\beta} > 1/2$ whenever $\alpha \not = \beta$. Let $y_\alpha = e \otimes x_\alpha$. Observe that $\norm{y_\alpha - y_\beta} > 1/2$ whenever $\alpha \not = \beta$ and that
    \begin{equation*}
        \norm{(U \hat{\otimes} I) (y_\alpha - y_\beta)} = \norm{(Ue) \otimes (x_\alpha - x_\beta)} = \norm{Ue} \norm{x_\alpha - x_\beta} \leq 2 \norm{Ue}.
    \end{equation*}
    Since $S$ has separable range, we can find $\alpha <  \beta < \omega_1$ such that $\norm{Sy_\alpha - Sy_\beta} < 1/(4\norm{T^{-1}})$. We have that
    \begin{align*}
        1/(2\norm{T^{-1}}) &\leq \norm{T(y_\alpha - y_\beta)} = \norm{(U \hat{\otimes} I) (y_\alpha - y_\beta) + S(y_\alpha - y_\beta)} \\
        &\leq \norm{(U \hat{\otimes} I) (y_\alpha - y_\beta)} + \norm{S(y_\alpha - y_\beta)} \leq 2 \norm{Ue} + 1/(4\norm{T^{-1}}),
    \end{align*}
    in other words $\norm{Ue} \geq 1/(8\norm{T^{-1}})$, so that $U$ is bounded below.

    We show now that $U[E]$ is dense in $F$. Thus, fix $f \in F$ and $\varepsilon > 0$, we will show that $d(f, U[E]) < \varepsilon$.
    
    As in the proof of Theorem \ref{th: main-tensor-prod} we can isometrically embed $F \hat{\otimes}_\varepsilon X$ into $C([0,1], X)$. Let $E_0$ be a countable dense set in $E$ and apply the first part of Lemma \ref{lmm: approx-coordinatewise} to obtain $x \in B_X$ and $\varphi \in B_{X^*}$, where we identify $S$ with $(j \hat{\otimes} I) \circ S$ whenever appropriate.

    Since $T$ is an isomorphism and $E_0$ is dense in $E$, we can find $e_1, \dots, e_n \in E_0$ and $x_1, \dots, x_n \in X$ such that
    \begin{equation*}
        \norm{(f \otimes x) - T \left(\sum_{i=1}^n  e_i \otimes x_i \right)} < \varepsilon/2.
    \end{equation*}
    Since $j \hat{\otimes} I$ is an isometry and $M_0$ is dense, this norm can be computed as
    \begin{align*}
        \sup_{k \in M_0} &\norm{\pi_k \left((j \hat{\otimes} I)(f \otimes x) - (j \hat{\otimes} I)\circ T \left(\sum_{i=1}^n  (e_i \otimes x_i) \right) \right)} \\
        &= \sup_{k \in M_0}\norm{j(f)(k)x - \sum_{i=1}^n  j(Ue_i)(k) x_i + \pi_k \circ (j \hat{\otimes} I) \circ S( (e_i \otimes x_i))} \\
        &\geq \frac{1}{2} \sup_{k \in M_0} |j(f)(k) - \sum_{i=1}^n \frac{\varphi(x_i)}{\varphi(x)}j(Ue_i)(k)| \\ 
        &\geq \frac{1}{2} d(j(f), j(U[E])) = \frac{1}{2}d(f, U[E]),
    \end{align*}
    which finishes the proof.
\end{proof}

%% file: 4-Applications.tex
\section{Applications} \label{sec: proof-applications}

We denote by $[\mathbb{N}]^{<\omega}$ and $[\mathbb{N}]^\omega$ the family of all finite and infinite subsets of $\mathbb{N}$, respectively. A family $\mathcal{A} \subseteq [\mathbb{N}]^\omega$ is called \emph{almost disjoint} if for any distinct $A, B \in \mathcal{A}$ then $A \cap B \in [\mathbb{N}]^{<\omega}$. In this section, we adopt the notation from \cite{koszmider2021banach} and write $K_\mathcal{A}$ for the locally compact Hausdorff space associated to the almost disjoint family $\mathcal{A} \subseteq [\mathbb{N}]^\omega$, as defined in \cite[Definition 3]{koszmider2021banach}. For convenience, we recall this definition below.

\begin{definition}
    Let $\mathcal{A} \subseteq [\mathbb{N}]^\omega$ be an almost disjoint family. Then $K_{\mathcal{A}}$ denotes the topological space consisting of distinct points $\{x_n : n \in \mathbb{N} \} \cup \{ y_A : A \in \mathcal{A} \}$, where $x_n$ is isolated for every $n \in \mathbb{N}$, and the sets 
    \begin{equation*}
        U(A,G) = \{ x_n : n \in A \setminus G \} \cup \{ y_A \}
    \end{equation*}  
    for $G \in [\mathbb{N}]^{<\omega}$ form a neighbourhood basis at each point $y_A$ for $A \in \mathcal{A}$. We write $U(A)$ for $U(A, \emptyset)$.  
\end{definition}

We begin by summarizing some fundamental properties of these spaces (see \cite[Lemma 4]{koszmider2021banach}) and extending certain results from Section 3.2 of the same paper for $C_0(K_\mathcal{A})$ to their natural counterparts in the context of $C_0(M \times K_\mathcal{A})$, where we recall that $M$ is a compact metric space.

\begin{lemma}\label{lmm: alm-disj-fam-spaces-recall}
    Let $\mathcal{A} \subseteq [\mathbb{N}]^\omega$ be an almost disjoint family. Then:
    \begin{itemize}
        \item $K_\mathcal{A}$ is a locally compact, scattered Hausdorff space.
        \item $K_\mathcal{A}$ is compact if and only if $\mathcal{A}$ and $\mathbb{N} \setminus \bigcup \mathcal{A}$ are both finite.
        \item $\{x_n : n \in \mathbb{N}\}$ is the set of isolated points of $K_\mathcal{A}$; it is dense in $K_\mathcal{A}$, and so $K_\mathcal{A}$ is separable. Hence $K_\mathcal{A}$ is metrizable if and only if it is second countable, if and only if $\mathcal{A}$ is countable.
        \item The subspace $K_\mathcal{A} \setminus \{x_n : n \in \mathbb{N}\} = \{y_A : A \in \mathcal{A}\}$ is closed and discrete.
        \item The sequence $(x_n)_{n \in A}$ converges to $y_A$ in $K_\mathcal{A}$ for every $A \in \mathcal{A}$.
 \end{itemize}
\end{lemma} 

\begin{definition}
    Let $\mathcal{A} \subseteq [\mathbb{N}]^\omega$ be an almost disjoint family. For $f \in C_0(M \times K_\mathcal{A})$ and $\mathcal{X} \subseteq C_0(M \times K_\mathcal{A})$, we define
    \begin{equation*}
        s(f) = \{A \in \mathcal{A}: \exists k \in M, f(k, y_A) \not = 0\}, \hspace{5pt} s(\mathcal{X}) = \bigcup \{s(f): f \in \mathcal{X}\}.
    \end{equation*}
\end{definition}

\begin{lemma}\label{lnmm: sep-supp-functions}
    The set $s(f)$ is countable for every $f \in C_0(M \times K_\mathcal{A})$.
\end{lemma}
\begin{proof}
    Let $M_0 \subseteq M$ be a countable dense set. We claim that
    \begin{equation*}
        s(f) = \{A \in \mathcal{A}: \exists k \in M_0, f(k, y_A) \not = 0\} = \bigcup_{k \in M_0} \{A \in \mathcal{A}: f(k, y_A) \not = 0\}.
    \end{equation*}
    Indeed, if $A \in s(f)$, then there exists $k \in M$ such that $f(k, y_A) \not = 0$, and taking $(k_n) \subseteq M_0$, $k_n \to k$ we have $f(k_n, y_A) \to f(k, y_A)$. So $f(k_N, y_A) \not = 0$ for some $N \in \mathbb{N}$

   Thus, it is enough to show that for each $k \in M_0$ the set $\{A \in \mathcal{A}: f(k, y_A) \not = 0\}$ is countable. Fix $k \in M_0$, and note that only finite subsets of the discrete closed set $\{(k, y_A): A \in \mathcal{A}\}$ are compact, while for each $\varepsilon > 0$, the set $\{(k, y_A) \in M \times K_\mathcal{A}: |f(k, y_A)| \geq \varepsilon \}$ is compact. The result follows.
\end{proof}

\begin{lemma}\label{lmm: separable-if-countable}
   Let $\mathcal{A} \subseteq [\mathbb{N}]^\omega$ be an almost disjoint family. A closed subspace $\mathcal{X}$ of $C_0(M \times K_\mathcal{A})$ is separable if and only if $s(\mathcal{X})$ is countable.
\end{lemma}

\begin{proof}
    Suppose that $\mathcal{X}$ is separable and let $\mathcal{X}_0$ be a countable dense subset of $\mathcal{X}$, then by Lemma \ref{lnmm: sep-supp-functions} $s(\mathcal{X}_0)$ is countable. It is easy to see that $s(\mathcal{X}) = s(\mathcal{X}_0)$, which proves the claim.

    On the other hand, if $s(\mathcal{X})$ is countable, then $\mathcal{X}$ is isomorphic to a subspace of
    \begin{equation*}
        C_0(M \times K_{s(\mathcal{X})}) \cong C(M) \hat{\otimes}_\varepsilon C_0(K_{s(\mathcal{X})}),
    \end{equation*}
    and by Lemma \ref{lmm: alm-disj-fam-spaces-recall}, $C_0(K_{s(\mathcal{X})})$ is separable and thus so is $C_0(M \times K_{s(\mathcal{X})})$.
\end{proof}

As usual, $K'$ denotes the Cantor--Bendixson derivative of a topological space $K$, that is, the set of accumulation points of $K$ equipped with the subspace topology. For any ordinal $\alpha$, we define
\begin{equation*}
    K^{(\alpha+1)} = \bigl(K^{(\alpha)}\bigr)', \quad \text{and if } \alpha \text{ is a limit ordinal, } K^{(\alpha)} = \bigcap_{\beta < \alpha} K^{(\beta)}.
\end{equation*}
The \emph{Cantor--Bendixson index} of $K$ is the smallest ordinal $\alpha$ such that $K^{(\alpha)} = \emptyset$; if no such ordinal exists, we say that the Cantor--Bendixson index is $\infty$. We observe that $K_\mathcal{A}$ has finite Cantor--Bendixson index for any almost disjoint family $\mathcal{A} \subseteq [\mathbb{N}]^\omega$.

\begin{lemma}\label{lmm: separable-contained-in-spaces}
    Suppose that $\mathcal{A} \subseteq [\mathbb{N}]^\omega$ is an almost disjoint family and $\mathcal{X} \subseteq C_0(M \times K_\mathcal{A})$ is separable. Then there is a closed subspace $\mathcal{Y}$ of $C_0(M \times K_\mathcal{A})$ such that $\mathcal{X} \subseteq \mathcal{Y}$ and $\mathcal{Y} \sim C_0( M \times \omega)$. 
    Specifically, when $M$ is infinite, $\mathcal{Y} \sim C(M)$, while for finite $M$, $\mathcal{Y} \sim c_0$.
\end{lemma}
\begin{proof}
    This follows by examining the second part of the proof of Lemma~\ref{lmm: separable-if-countable} and observing that, since $K_{s(\mathcal{X})}$ has finite Cantor--Bendixson index, the classical classification of continuous functions on countable compact Hausdorff spaces \cite{bessaga1960spaces} implies that $C_0(K_{s(\mathcal{X})}) \sim C_0(\omega)$, or, if $s(\mathcal{X})$ is finite, $C_0(K_{s(\mathcal{X})}) \sim \mathbb{K}^{|s(\mathcal{X})|} \subseteq C_0(\omega)$.
\end{proof}

Our next lemma is a reformulation of \cite[Proposition 3.11]{candido2021banach} in the context of spaces of the form  $C_0(M \times K_\mathcal{A})$. Since the proof is essentially identical, we omit it here.

\begin{lemma}\label{lmm: projections-separable}
    Let $P: C_0(M \times K_\mathcal{A}) \to C_0(M \times K_\mathcal{A})$ be a projection and suppose $P = Q \hat{\otimes} I + S$ for an operator $Q: C(M) \to C(M)$ and an operator $S$ with separable range. Then $Q$ is also a projection.
\end{lemma}

Lastly, we will need the following well-known result, see for example \cite[Lemma 3.6 (ii)]{laustsen2002maximal}. Recall that an operator $P: \mathcal{X} \to \mathcal{X}$ is said to \emph{factor through} $\mathcal{Y}$ if there exist operators $Q: \mathcal{X} \to \mathcal{Y}$ and $R: \mathcal{Y} \to \mathcal{X}$ such that $P = RQ$.

\begin{lemma}\label{lmm: factorization-projections}
    Let $\mathcal{X}, \mathcal{Y}$ be Banach spaces and $P: \mathcal{X} \to \mathcal{X}$ be a projection. Suppose $P$ factors through $\mathcal{Y}$, then $P[\mathcal{X}]$ is isomorphic to a complemented subspace of $\mathcal{Y}$.
\end{lemma}

\begin{proposition}\label{prop: fin-dimensional-case}
    Let $n \in \mathbb{N}$ and suppose that $C_0(K_{\mathcal{A}})$ has few operators. Then every infinite-dimensional complemented subspace of $C_0(K_\mathcal{A})^n$ is isomorphic to either $c_0$ or $C_0(K_\mathcal{A})^m$ for some $m \leq n$.
\end{proposition}
\begin{proof}
    We proceed by induction, the case $n=1$ follows by \cite[Lemma 4]{koszmider2005decompositions}. Thus, assume the result is true for $m < n$, and we will show it holds for $n$.

    Let $P: C_0(K_\mathcal{A})^n \to C_0(K_\mathcal{A})^n$ be a projection with infinite-dimensional image, and write $P = Q \hat{\otimes} I - S$ using Corollary \ref{corollary-operators-CK}, where $S$ has separable range. By Lemma \ref{lmm: projections-separable}, the operator $Q: \ell_\infty^n \to \ell_\infty^n$ is also a projection and let $d = \dim Q[\ell_\infty^n]$.  
    Since $S$ factors through $c_0$ by Lemma \ref{lmm: separable-contained-in-spaces}, it follows that $P$ factors through $(Q[\ell_\infty^n] \hat{\otimes}_{\varepsilon} C_0(K_\mathcal{A}))\oplus c_0 \sim C_0(K_\mathcal{A})^{d} \oplus c_0$. We distinguish three cases: $d = 0$, $0 < d < n$ or $d = n$.  
    
    In the first case, $Q = 0$, so $P = S$ which factors through $c_0$, so that $P[C_0(K_\mathcal{A})^n]$ is isomorphic to a complemented subspace of $c_0$ and thus isomorphic to $c_0$ \cite[Theorem 1]{pelczynski1960projections}.
    
    If $0 < d < n$, $P$ factors through $C_0(K_\mathcal{A})^{d} \oplus c_0 \sim C_0(K_\mathcal{A})^{d}$, so that an application of the induction hypothesis, together with Lemma \ref{lmm: factorization-projections}, gives the result.  
    
    In the case $d = n$ then $Q = I$, so that $P = I - S$. Observe that $I - P = S$, implying that $S$ is a projection that factors through $c_0$. By Lemma \ref{lmm: factorization-projections}, the space $S[C_0(K_\mathcal{A})^n]$ is isomorphic to a complemented subspace of $c_0$ and thus either finite-dimensional or isomorphic to $c_0$.
    
    Since $P[C_0(K_\mathcal{A})^n]$ is an infinite-dimensional complemented subspace of a space of continuous functions over a scattered locally compact space, it contains a complemented copy of $c_0$ \cite{lotz1979semi}. Therefore
    \begin{equation*}
        C_0(K_\mathcal{A})^n \sim P[C_0(K_\mathcal{A})^n] \oplus S[C_0(K_\mathcal{A})^n] \sim P[C_0(K_\mathcal{A})^n],
    \end{equation*}
    which completes the proof.
\end{proof}

We are ready to provide the classification of the complemented subspaces of $C_0(\omega \times K_\mathcal{A})$.

\begin{proof}[Proof of Theorem \ref{th: complemented-subspaces}]
     Let $P: C_0(\omega \times K_\mathcal{A}) \to C_0(\omega \times K_\mathcal{A})$ be a projection and express $P = Q \hat{\otimes} I - S$ according to Corollary \ref{corollary-operators-CK}, so that by Lemma \ref{lmm: projections-separable} $Q$ is also a projection. If $Q[C_0(\omega)]$ is finite dimensional, then $P$ factors through $C_0(K_\mathcal{A})^{\dim Q[C_0(\omega)]+1}$, and the result follows from Lemma \ref{lmm: factorization-projections} and Proposition \ref{prop: fin-dimensional-case}.

    Otherwise, $Q[C_0(\omega)] \sim C_0(\omega)$, so after conjugation if necessary, we can assume $Q = I$, which gives $P = I - S$. Arguing as in the proof of Proposition \ref{prop: fin-dimensional-case} gives $P[C_0(\omega \times K_\mathcal{A})] \sim C_0(\omega \times K_\mathcal{A})$.
\end{proof}

Finally, we can classify the closed ideals of $\mathscr{B}(C_0(\omega \times K_\mathcal{A}))$.

\begin{proof}[Proof of Theorem \ref{th: ideal-lattice}]
Observe that $C_0(\omega \times K_\mathcal{A})$ has the bounded approximation property, being a $\mathscr{L}_\infty$-space, and thus $\mathscr{K}(C_0(\omega \times K_\mathcal{A}))$ is the minimal non-zero closed ideal. By Lemma \ref{lmm: separable-contained-in-spaces}, operators with separable range factor through $c_0$, so that $\mathscr{G}_{c_0}(C_0(\omega \times K_\mathcal{A})) = \mathscr{X}(C_0(\omega \times K_\mathcal{A}))$. In particular $\mathscr{G}_{c_0}(C_0(\omega \times K_\mathcal{A}))$ is a closed operator ideal. 

Since $\omega \times K_\mathcal{A}$ is scattered, it is well-known that non-compact operators in $C_0(\omega \times K_\mathcal{A})$ factor the identity of $c_0$ \cite[Proposition 5.4 (ii)]{kania2014ideal}, so that any ideal $\mathscr{I}$ which is not contained in the ideal of compact operators, also contains the ideal of separable operators. Therefore for any closed ideal $\mathscr{I}$ of $\mathscr{B}(C_0(\omega \times K_\mathcal{A}))$, $\mathscr{I} \not = \mathscr{K}(C_0(\omega \times K_\mathcal{A}))$, we have
\begin{equation*}
    \{0\} \subsetneq \mathscr{K}(C_0(\omega \times K_\mathcal{A})) \subsetneq \mathscr{X}(C_0(\omega \times K_\mathcal{A})) \subseteq \mathscr{I}.
\end{equation*}
Since $\varphi: \mathscr{B}(C_0(\omega \times K_\mathcal{A})) \to \mathscr{B}(c_0)$, defined by $\varphi(T) = \varphi(U \hat{\otimes} I + S) = U$, is a surjective algebra homomorphism with $\ker{\varphi} = \mathscr{X}((C_0(\omega \times K_\mathcal{A}))$, it follows that there is a one-to-one correspondence between the closed ideals of $\mathscr{B}(c_0)$ and the closed ideals of $\mathscr{B}(C_0(\omega \times K_\mathcal{A}))$ which contain the ideal $\mathscr{X}(C_0(\omega \times K_\mathcal{A}))$. The only closed operator ideals on $c_0$ are $\{0\} \subsetneq \mathscr{K}(c_0) \subsetneq \mathscr{B}(c_0)$ \cite{gohberg1967normally}; we obtain the rest of the ideals of $\mathscr{B}(C_0(\omega \times K_\mathcal{A}))$ as the preimages of those. The only non-trivial case is
\begin{align*}
   \varphi^{-1}(\mathscr{K}(c_0) ) &= \{Q\hat{\otimes} I+S: Q \in \mathscr{K}(c_0), S \in \mathscr{X}(C_0(\omega \times K_\mathcal{A})) \} \\
   &= \overline{\mathscr{G}}_{C_0(K_\mathcal{A})}(C_0(\omega \times K_\mathcal{A})),
\end{align*}
which follows since finite rank operators uniformly approximate every compact operator in $c_0$, while separable range operators on $C_0(\omega \times K_\mathcal{A})$ factor through $c_0$, and thus through $C_0(K_\mathcal{A})$.
\end{proof}